# A Census of Non-Reconstructable Digraphs, II:

# A Family of Tournaments

Paul K. Stockmeyer
Department of Computer Science
College of William and Mary
Williamsburg, Virginia 23185

ABSTRACT: Continuing the program begun in the first paper of this series, we present a pair of non-reconstructable tournaments on p points for each $p = 2^n \geq 4$.

## 1. Introduction

The graph reconstruction conjecture continues to be one of the most challenging open problems in graph theory, and the recent survey article by Manvel [2] offers little hope that it will be settled soon. Even though he lists over 150 papers published in just the last 10 years related to this conjecture, most of these concern very special cases, or else they explore analogous conjectures for a variety of other graph-like structures. The original conjecture for graphs appears to be as intractable as ever.

The situation for directed graphs is quite different. It has been known for several years that there are infinite families of digraphs that contradict the reconstruction conjecture. In the directed case, then, the interesting activities are to catalog and study the rich structure of the digraphs that satisfy the following definitions.

DEFINITION 1: Two digraphs G and G* are called <u>hypomorphic</u> if their points can be labeled $\{v_1, v_2, ..., v_p\}$ and $\{u_1, u_2, ..., u_p\}$, respectively, so that the point-deleted subgraphs $G - v_i$ and $G^* - u_i$ are isomorphic for each i from 1 to p. The digraphs are called <u>non-reconstructable</u> if they are hypomorphic but not isomorphic.

In [3] we presented six pairs of non-reconstructable digraphs on p points for each p of the form $p = 2^m + 2^n$ with $0 \leq n < m$, including one pair of tournaments. We also noted that three of these pairs remain non-reconstructable in the case $p = 2^n$. The two tournaments, though, happen to be isomorphic in this case, and hence reconstructable. In this paper we remedy the situation by presenting a new family of non-reconstructable tournaments, one pair for each order $p = 2^n \geq 4$.

## 2. The Matrices

We begin by defining two p × p matrices $M_p$ and $M_p^*$ for each $p = 2^n \geq 4$. For the moment, these matrices can be viewed as the adjacency matrices of two edge-weighted digraphs. Later we will see how to convert these to the adjacency matrices of non-reconstructable tournaments.

DEFINITION 2: Matrices $M_4$ and $M_4^*$ are by definition the matrices displayed in Figure 1. For $p = 2^n \geq 8$, $M_p$ and $M_p^*$ each consist of a p/4 × p/4 array of blocks of size 4 × 4, defined by

$$\text{block } i,j \text{ of } M_p = \begin{cases} M_4 & \text{if } j-i = 0 \\ -M_4 + 4I & \text{if } j-i \equiv 1 \pmod{4} \\ -M_4 - 4I & \text{if } j-i \equiv -1 \pmod{4} \\ M_4 + (x+4)I & \text{if } j-i = y \cdot 2^x \text{ with } x \geq 1 \text{ and } y \equiv 1 \pmod{4} \\ M_4 - (x+4)I & \text{if } j-i = y \cdot 2^x \text{ with } x \geq 1 \text{ and } y \equiv -1 \pmod{4} \end{cases}$$

and

$$\text{block } i,j \text{ of } M_p^* = \begin{cases} M_4^* & \text{if } j-i = 0 \\ -M_4^* - 4I & \text{if } j-i \equiv 1 \pmod{4} \\ -M_4^* + 4I & \text{if } j-i \equiv -1 \pmod{4} \\ M_4^* - (x+4)I & \text{if } j-i = y \cdot 2^x \text{ with } x \geq 1 \text{ and } y \equiv 1 \pmod{4} \\ M_4^* + (x+4)I & \text{if } j-i = y \cdot 2^x \text{ with } x \geq 1 \text{ and } y \equiv -1 \pmod{4} \end{cases}$$

where $I$ is the $4 \times 4$ identity matrix.

We illustrate this definition by displaying the $8 \times 8$ and the $16 \times 16$ matrices in Figures 2 and 3 respectively.

We next list some properties of these matrices that will be needed in later sections. These results follow easily from the definitions, and the proofs are left to the reader.

LEMMA 1: (a). For all $p = 2^n \geq 8$, both the upper-left and the lower-right $p/2 \times p/2$ submatrices of $M_p$ are copies of $M_{p/2}$. The same is true for $M_p^*$ and $M_{p/2}^*$. More formally,

$$M_{p/2}[i,j] = M_p[i,j] = M_p[i + p/2, j + p/2]$$

and

$$M_{p/2}^*[i,j] = M_p^*[i,j] = M_p^*[i + p/2, j + p/2]$$

for all $i$ and $j$ from 1 to $p/2$.

(b). For all distinct i and j from 1 to 4,

$$M_8[i,j] = -M_8[i, j + 4] = -M_8[i + 4, j]$$

and

$$M_8^*[i,j] = -M_8^*[i, j + 4] = -M_8^*[i + 4, j].$$

(c). For all $p = 2^n \geq 16$, and all distinct i and j from 1 to p/2,

$$M_p[i,j] = \pm M_p[i, j + p/2] = \pm M_p[i + p/2, j]$$

and

$$M_p^*[i,j] = \pm M_p^*[i, j + p/2] = \pm M_p^*[i + p/2, j],$$

with the + sign used except when $|j-i| = p/4$.

(d). For all $p = 2^n \geq 8$ and all i from 1 to p/2,

$$M_p[i, i + p/2] = M_p^*[i + p/2, i] = (n + 1)$$

and

$$M_p[i + p/2, i] = M_p^*[i, i + p/2] = -(n + 1).$$

## 3. The Mappings

Our goal in the next section is to prove that the edge-weighted digraphs represented by the matrices $M_p$ and $M_p^*$ are hypomorphic. Here we present the isomorphisms that will accomplish this goal. Specifically, for each $p = 2^n \geq 4$ and each k from 1 to p, we define a bijection $\sigma_{p,k}$ from the set $\{1, 2, \ldots, p\} - \{k\}$ onto itself that will later serve as the k-th point-deleted isomorphism. The definition is inductive on p.

DEFINITION 3: The mappings $\sigma_{4,k}$, for k from 1 to 4, are defined by their values listed in Figure 4. For $p = 2^n \geq 8$, $\sigma_{p,k}$ is defined by

$$\sigma_{p,k}(i) = \begin{cases} \sigma_{p/2,k}(i) + p/2 & \text{if } k \leq p/2 \text{ and } i \leq p/2 \\ & \text{(with } i \neq k\text{)} \\ \sigma_{p/2,k}(i - p/2) & \text{if } k \leq p/2 \text{ and } i > p/2 \\ & \text{but } i \neq k + p/2 \\ i & \text{if } k \leq p/2 \text{ and } i = k + p/2 \\ \sigma_{p/2,k-p/2}(i) + p/2 & \text{if } k > p/2 \text{ and } i \leq p/2 \\ & \text{but } i \neq k - p/2 \\ \sigma_{p/2,k-p/2}(i - p/2) & \text{if } k > p/2 \text{ and } i > p/2 \\ & \text{(with } i \neq k\text{)} \\ i & \text{if } k > p/2 \text{ and } i = k - p/2. \end{cases}$$

Note that $\sigma_{p,k}(i)$ is undefined for $i = k$.

For convenience, the values for the mappings $\sigma_{8,k}$ and $\sigma_{16,k}$ are given in Figures 5 and 6.

Some of the properties of these mappings that will be needed later are listed in the following lemma. Parts (a), (b), and (c) are immediate consequences of Definition 3, while part (d) follows easily from part (b) and the bijective property of the $\sigma$ functions.

LEMMA 2: (a). For all $p = 2^n \geq 8$, all k from 1 to p/2, and all i from 1 to p distinct from k and k + p/2,

$$\sigma_{p,k}(i) = \sigma_{p,k+p/2}(i).$$

(b). For all $p = 2^n \geq 8$, all i from 1 to p/2, and all k from 1 to p distinct from i and i + p/2,

$$\sigma_{p,k}(i + p/2) = \sigma_{p,k}(i) - p/2.$$

(c). For all $p = 2^n \geq 8$, and for all distinct i and j from 1 to p,

$j = i \pm p/2$ if and only if $\sigma_{p,i}(j) = i \pm p/2$,

with the same sign in both cases.

(d). For all $p = 2^n \geq 8$, all k from 1 to p, and all i and j from 1 to p distinct from k,

$j = i \pm p/2$ if and only if $\sigma_{p,k}(i) = \sigma_{p,k}(j) \pm p/2$,

again with the same sign in both cases.

## 4. The Main Results

Before proving that our edge-weighted digraphs are hypomorphic, we first present a lemma that will significantly simplify that proof in some cases.

LEMMA 3: (a). For all distinct i and j from 1 to 4,

$$M_4[i,j] = -M_4^*[i,\sigma_{4,i}(j)] = -M_4^*[\sigma_{4,j}(i),j].$$

(b). For all $p = 2^n \geq 8$, and all distinct i and j from 1 to p,

$$M_p[i,j] = \pm M_p^*[i,\sigma_{p,i}(j)] = \pm M_p^*[\sigma_{p,j}(i),j],$$

with the + sign to be used except when $|j-i| = p/2$.

PROOF: Part (a) can be confirmed easily by inspection, as can part (b) for $p = 8$. We shall proceed by induction, verifying the claim in part (b) for an arbitrary matrix size $p \geq 16$ under the assumption that it is true for matrix size $p/2$. Unfortunately, it seems necessary to examine various cases, depending on the relative sizes of i and j. We will prove only the first equality; the proof of the second is essentially the same. In most cases, the derivation makes use of Lemma 1 (a) twice, Lemma 1 (c), the induction hypothesis, and the definition of the σ functions, in some order.

Case A: Both i and j are at most $p/2$. In this case we have

$$M_p[i,j] = M_{p/2}[i,j]$$
$$= \pm M^*_{p/2}[i, \sigma_{p/2,i}(j)]$$
$$= \pm M^*_p[i, \sigma_{p/2,i}(j)]$$
$$= M^*_p[i, \sigma_{p/2,i}(j) + p/2]$$
$$= M^*_p[i, \sigma_{p,i}(j)].$$

The induction hypothesis dictates that the negative sign must be introduced if and only if $|j-i| = p/4$. A later sign change takes place, according to Lemma 1 (c), precisely when $|\sigma_{p/2,i}(j) - i| = p/4$. Lemma 2 (c) assures us that these two conditions are equivalent, so that the net effect is no sign change.

Case B: Both i and j are greater than $p/2$. This case is similar to Case A. We have

$$M_p[i,j] = M_{p/2}[i - p/2, j - p/2]$$
$$= \pm M^*_{p/2}[i - p/2, \sigma_{p/2,i-p/2}(j - p/2)]$$
$$= \pm M^*_p[i - p/2, \sigma_{p/2,i-p/2}(j - p/2)]$$
$$= M^*_p[i, \sigma_{p/2,i-p/2}(j - p/2)]$$
$$= M^*_p[i, \sigma_{p,i}(j)].$$

The first sign change is made if $|(j-p/2) - (i-p/2)| = p/4$, and the sign changes back if $|\sigma_{p/2,i-p/2}(j-p/2) - (i-p/2)| = p/4$. As before, Lemma 2 (c) assures us that these conditions are equivalent.

Case C: Only i is at most $p/2$; j is greater than $p/2$. This case requires further subdivision as follows:

Case C.1: $j \neq i + p/2$. Here we have

$$\begin{aligned} M_p[i,j] &= \pm M_p[i, j - p/2] \\ &= \pm M_{p/2}[i, j - p/2] \\ &= M^*_{p/2}[i, \sigma_{p/2,i}(j - p/2)] \\ &= M^*_p[i, \sigma_{p/2,i}(j - p/2)] \\ &= M^*_p[i, \sigma_{p,i}(j)]. \end{aligned}$$

The negative sign appears and then vanishes in the case $|(j-p/2) - i| = p/4$.

Case C.2: $j = i + p/2$. In this case Lemma 1 (d) and Definition 3 yield

$$\begin{aligned} M_p[i,j] &= (n + 1) \\ &= -M^*_p[i,j] \\ &= -M^*_p[i, \sigma_{p,i}(j)]. \end{aligned}$$

Case D: Only $j$ is at most $p/2$; $i$ is greater than $p/2$. This is similar to case C, and again we have two sub-cases.

Case D.1: $i \neq j + p/2$. Then

$$\begin{aligned} M_p[i,j] &= \pm M_p[i - p/2, j] \\ &= \pm M_{p/2}[i - p/2, j] \\ &= M^*_{p/2}[i - p/2, \sigma_{p/2, i-p/2}(j)] \\ &= M^*_p[i, \sigma_{p/2, i-p/2}(j) + p/2] \\ &= M^*_p[i, \sigma_{p,i}(j)], \end{aligned}$$

using the negative signs if $|j - (i-p/2)| = p/4$.

Case D.2: $i = j + p/2$. Again, Lemma 1 (d) and Definition 3 yield

$$\begin{aligned} M_p[i,j] &= -(n + 1) \\ &= -M^*_p[i,j] \\ &= -M^*_p[i, \sigma_{p,i}(j)]. \end{aligned}$$

We are now ready to prove that our edge-weighted digraphs are hypomorphic.

THEOREM 1: For each $p = 2^n \geq 4$, the edge-weighted digraphs with adjacency matrices $M_p$ and $M_p^*$ are hypomorphic under the set of mappings $\sigma_{p,k}$. More specifically, for each k from 1 to p, and for all i and j from 1 to p distinct from k,

$$M_p[i,j] = M_p^*[\sigma_{p,k}(i), \sigma_{p,k}(j)].$$

PROOF: This proof also is by induction on p. The case p = 4 can be confirmed easily by inspection. We shall verify the truth of the claim for an arbitrary matrix size $p \geq 8$, under the assumption that it is true for matrix size p/2. Again, it seems necessary to examine various cases, depending on the relative sizes of i, j, and k. The claim is obviously true when i = j, since the diagonal elements of both matrices are all 0. Thus we shall assume in what follows that i is distinct from j. Also, we can restrict our attention to the cases where $k \leq p/2$; the proofs for $k > p/2$ are similar.

Case A: Both i and j are at most p/2. Using Lemma 1 (a), the induction hypothesis, and the definition of $\sigma_{p,k}$, we have

$$\begin{aligned}
M_p[i,j] &= M_{p/2}[i,j] \\
&= M_{p/2}^*[\sigma_{p/2,k}(i), \sigma_{p/2,k}(j)] \\
&= M_p^*[\sigma_{p/2,k}(i) + p/2, \sigma_{p/2,k}(j) + p/2] \\
&= M_p^*[\sigma_{p,k}(i), \sigma_{p,k}(j)].
\end{aligned}$$

Case B: Both i and j are greater than p/2. This case requires further subdivision because of the special definition of $\sigma_{p,k}(i)$ when $i = k + p/2$.

Case B.1: Neither i nor j equals $k + p/2$. Essentially the same reasoning as in Case A yields

$$\begin{aligned}
M_p[i,j] &= M_{p/2}[i - p/2, j - p/2] \\
&= M^*_{p/2}[\sigma_{p/2,k}(i - p/2), \sigma_{p/2,k}(i - p/2)] \\
&= M^*_p[\sigma_{p/2,k}(i - p/2), \sigma_{p/2,k}(j - p/2)] \\
&= M^*_p[\sigma_{p,k}(i), \sigma_{p,k}(j)].
\end{aligned}$$

Case B.2: $i = k + p/2$. (We omit the proof of the case $j = k + p/2$, which is essentially the same). This is where we use Lemma 3 and Lemma 2 (a).

$$\begin{aligned}
M_p[i,j] &= M^*_p[i, \sigma_{p,i}(j)] \\
&= M^*_p[i, \sigma_{p,k}(j)] \\
&= M^*_p[\sigma_{p,k}(i), \sigma_{p,k}(j)].
\end{aligned}$$

Case C: Only i is at most $p/2$; j is greater than $p/2$. Here we have three subcases.

Case C.1: j equals neither $i + p/2$ nor $k + p/2$. Using Lemma 1 (a) and 1 (b) (if $p = 8$) or 1 (c) (if $p \geq 16$), we have

$$\begin{aligned}
M_p[i,j] &= \pm M_p[i, j - p/2] \\
&= \pm M_{p/2}[i, j - p/2] \\
&= \pm M^*_{p/2}[\sigma_{p/2,k}(i), \sigma_{p/2,k}(j - p/2)] \\
&= \pm M^*_p[\sigma_{p/2,k}(i), \sigma_{p/2,k}(j - p/2)] \\
&= M^*_p[\sigma_{p/2,k}(i) + p/2, \sigma_{p/2,k}(j - p/2)] \\
&= M^*_p[\sigma_{p,k}(i), \sigma_{p,k}(j)].
\end{aligned}$$

Here the first sign change takes place if $p = 8$ or if $|(j-p/2) - i| = p/4$ and $p \geq 16$. The sign changes back if $p = 8$ or if $|\sigma_{p/2,k}(j-p/2) - \sigma_{p/2,k}(i)| = p/4$ with $p \geq 16$. This time it is Lemma 2 (d) that asserts the equivalence of these conditions.

Case C.2: $j = k + p/2$. This is similar to Case B.2. Using Lemma 3 again we have

$$M_p[i,j] = M_p^*[\sigma_{p,j}(i),j]$$
$$= M_p^*[\sigma_{p,k}(i),j]$$
$$= M_p^*[\sigma_{p,k}(i),\sigma_{p,k}(j)].$$

Case C.3: $j = i + p/2$. Using Lemmas 1 (d) and 2 (d) we have

$$M_p[i,j] = (n + 1)$$
$$= M_p^*[\sigma_{p,k}(i),\sigma_{p,k}(j)].$$

The final Case D, with i greater than p/2 and j at most p/2, is not essentially different from Case C. The proof is left for the reader to construct.

It will follow from the results of the next section that the edge-weighted digraphs represented by $M_p$ and $M_p^*$ are not isomorphic. We thus conclude that they form an infinite family of non-reconstructable edge-weighted digraphs.

## 5. Forming Digraphs

While non-reconstructable edge-weighted digraphs are interesting, we are primarily interested in standard digraphs. In this section we will show how to convert $M_p$ and $M_p^*$ into the adjacency matrices of non-reconstructable digraphs.

We can certainly convert $M_p$ and $M_p^*$ into the adjacency matrices of some digraphs by replacing each entry with either 0 or 1. Moreover, if for each number from $-(n+1)$ to $(n+1)$, we assign the same binary value to all occurrences of that number, the resulting digraphs will clearly be hypomorphic, as above. We call such an assignment <u>proper</u>. For example, suppose we replace all positive entries in $M_4$ and $M_4^*$ with 1, and all negative entries with 0. The matrices resulting from this proper assignment represent the 4-point non-reconstructable tournaments first displayed in [1]. The tricky part is to choose a proper assignment for which the resulting digraphs are non-isomorphic.

Consider the mapping $\sigma_{p,1}$ for some $p = 2^n \geq 8$. We can extend this function to a bijection on $\{1, 2, \ldots, p\}$ by setting $\sigma_{p,1}(1) = 1$. Then from Lemma 3 and Theorem 1, we know that

$$M_p[i,j] = M_p^*[\sigma_{p,1}(i), \sigma_{p,1}(j)]$$

for all i and j from 1 to p <u>except</u> for the cases

$$M_p[1, 1 + p/2] = (n + 1) = -M_p^*[\sigma_{p,1}(1), \sigma_{p,1}(1 + p/2)]$$

and

$$M_p[1 + p/2, 1] = -(n + 1) = -M_p^*[\sigma_{p,1}(1 + p/2), \sigma_{p,1}(1)].$$

Thus any proper assignment of binary values to the entries of $M_p$ and $M_p^*$ that assigns 0 to both (n+1) and -(n+1), or that assigns 1 to both (n+1) and -(n+1), will result in isomorphic digraphs, joined by the extension of the mapping $\sigma_{p,1}$. It is necessary to assign different binary values to (n+1) and -(n+1) in order to produce non-isomorphic digraphs.

Now consider the permutation $(1, 1 + p/2)(2, 2 + p/2) \ldots (p/2, p)$, again with $p = 2^n \geq 8$. When this permutation is applied to the rows and columns of $M_p$ and $M_p^*$, the only effect is to interchange the entries (n+1) and -(n+1) wherever they occur. Thus any proper assignment of binary values to the entries of $M_p$ and $M_p^*$ that replaces (n+1) with 1 and -(n+1) with 0 will result in the same digraph pair as the assignment that differs from this one only in that it replaces (n+1) with 0 and -(n+1) with 1. Combining this observation with the one in the previous paragraph, we conclude that in constructing non-isomorphic digraphs from the matrices $M_p$ and $M_p^*$, we may assume without loss of generality that all entries of (n+1) are replaced with 1 and that all entries of -(n+1) are replaced with 0.

Although the details are rather tedious, one can prove similar claims about the entries ±4, ±5, ... ±n. We might as well go ahead and replace all entries greater than or equal to 4 with 1, and all entries less than or equal to -4 with 0. Any non-isomorphic digraph pair that can be created from a proper assignment of binary values to the entries of $M_p$ and $M_p^*$ will arise in this way. Our only real choices are in assigning values to ±1, ±2, and ±3.

Without further ado, we now present our family of non-reconstructable tournaments.

THEOREM 2: For each $p = 2^n \geq 4$, let $G_p$ and $G_p^*$ be the tournaments with adjacency matrices obtained from $M_p$ and $M_p^*$, respectively, by replacing all positive entries with 1 and all negative entries with 0. Then $G_p$ and $G_p^*$ are non-isomorphic and hence non-reconstructable.

PROOF: Once again the proof is by induction on p. $G_4$ and $G_4^*$ are the well-known non-reconstructable tournaments mentioned above. For $p \geq 8$, we note that the first $p/2$ points of $G_p$ and the last $p/2$ points of $G_p^*$ have score (outdegree) of $p/2$, while the last $p/2$ points of $G_p$ and the first $p/2$ points of $G_p^*$ each have score $p/2 - 1$. Thus any isomorphism from $G_p$ to $G_p^*$ would have to map the first $p/2$ points of $G_p$ onto the last $p/2$ points of $G_p^*$. We know from Lemma 1 (a), though, that these point sets induce subtournaments isomorphic to $G_{p/2}$ and $G_{p/2}^*$, respectively. By hypothesis, the tournaments $G_{p/2}$ and $G_{p/2}^*$ are non-isomorphic. Hence $G_p$ and $G_p^*$ are non-isomorphic as well.

For completeness, we note that $M_p$ and $M_p^*$ can be used to produce a second pair of non-reconstructable digraphs. For this pair, assign the value 1 to entries 1, -1, and all entries $\geq 3$, and assign the value 0 to 2, -2, and all entries $\leq -3$. The proof of Theorem 2 can be modified easily to show that the resulting digraphs are non-reconstructable. This is not a new pair, however. These digraphs are in fact isomorphic to the digraphs $C_p$ anc $C_p^*$ described in [3] and generated there using a different set of point-deleted isomorphisms.

```
 0  1  2  3          0 -2 -3 -1
-1  0  3 -2          2  0  1 -3
-2 -3  0  1          3 -1  0  2
-3  2 -1  0          1  3 -2  0

   Matrix M₄          Matrix M₄*
```

Figure 1.  The smallest matrices.

```
 0  1  2  3  4 -1 -2 -3          0 -2 -3 -1 -4  2  3  1
-1  0  3 -2  1  4 -3  2          2  0  1 -3 -2 -4 -1  3
-2 -3  0  1  2  3  4 -1          3 -1  0  2 -3  1 -4 -2
-3  2 -1  0  3 -2  1  4          1  3 -2  0 -1 -3  2 -4
-4 -1 -2 -3  0  1  2  3          4  2  3  1  0 -2 -3 -1
 1 -4 -3  2 -1  0  3 -2         -2  4 -1  3  2  0  1 -3
 2  3 -4 -1 -2 -3  0  1         -3  1  4 -2  3 -1  0  2
 3 -2  1 -4 -3  2 -1  0         -1 -3  2  4  1  3 -2  0

           Matrix M₈                        Matrix M₈*
```

Figure 2.  The matrices of order 8.

```
 0  1  2  3  4 -1 -2 -3  5  1  2  3 -4 -1 -2 -3
-1  0  3 -2  1  4 -3  2 -1  5  3 -2  1 -4 -3  2
-2 -3  0  1  2  3  4 -1 -2 -3  5  1  2  3 -4 -1
-3  2 -1  0  3 -2  1  4 -3  2 -1  5  3 -2  1 -4
-4 -1 -2 -3  0  1  2  3  4 -1 -2 -3  5  1  2  3
 1 -4 -3  2 -1  0  3 -2  1  4 -3  2 -1  5  3 -2
 2  3 -4 -1 -2 -3  0  1  2  3  4 -1 -2 -3  5  1
 3 -2  1 -4 -3  2 -1  0  3 -2  1  4 -3  2 -1  5
-5  1  2  3 -4 -1 -2 -3  0  1  2  3  4 -1 -2 -3
-1 -5  3 -2  1 -4 -3  2 -1  0  3 -2  1  4 -3  2
-2 -3 -5  1  2  3 -4 -1 -2 -3  0  1  2  3  4 -1
-3  2 -1 -5  3 -2  1 -4 -3  2 -1  0  3 -2  1  4
 4 -1 -2 -3 -5  1  2  3 -4 -1 -2 -3  0  1  2  3
 1  4 -3  2 -1 -5  3 -2  1 -4 -3  2 -1  0  3 -2
 2  3  4 -1 -2 -3 -5  1  2  3 -4 -1 -2 -3  0  1
 3 -2  1  4 -3  2 -1 -5  3 -2  1 -4 -3  2 -1  0
```

Matrix $M_{16}$

```
 0 -2 -3 -1 -4  2  3  1 -5 -2 -3 -1  4  2  3  1
 2  0  1 -3 -2 -4 -1  3  2 -5  1 -3 -2  4 -1  3
 3 -1  0  2 -3  1 -4 -2  3 -1 -5  2 -3  1  4 -2
 1  3 -2  0 -1 -3  2 -4  1  3 -2 -5 -1 -3  2  4
 4  2  3  1  0 -2 -3 -1 -4  2  3  1 -5 -2 -3 -1
-2  4 -1  3  2  0  1 -3 -2 -4 -1  3  2 -5  1 -3
-3  1  4 -2  3 -1  0  2 -3  1 -4 -2  3 -1 -5  2
-1 -3  2  4  1  3 -2  0 -1 -3  2 -4  1  3 -2 -5
 5 -2 -3 -1  4  2  3  1  0 -2 -3 -1 -4  2  3  1
 2  5  1 -3 -2  4 -1  3  2  0  1 -3 -2 -4 -1  3
 3 -1  5  2 -3  1  4 -2  3 -1  0  2 -3  1 -4 -2
 1  3 -2  5 -1 -3  2  4  1  3 -2  0 -1 -3  2 -4
-4  2  3  1  5 -2 -3 -1  4  2  3  1  0 -2 -3 -1
-2 -4 -1  3  2  5  1 -3 -2  4 -1  3  2  0  1 -3
-3  1 -4 -2  3 -1  5  2 -3  1  4 -2  3 -1  0  2
-1 -3  2 -4  1  3 -2  5 -1 -3  2  4  1  3 -2  0
```

Matrix $M^*_{16}$

Figure 3. The matrices of order 16.

| i | $\sigma_{4,1}(i)$ | $\sigma_{4,2}(i)$ | $\sigma_{4,3}(i)$ | $\sigma_{4,4}(i)$ |
|---|---|---|---|---|
| 1 | X | 3 | 4 | 2 |
| 2 | 4 | X | 1 | 3 |
| 3 | 2 | 4 | X | 1 |
| 4 | 3 | 1 | 2 | X |

Figure 4. The mappings $\sigma_{4,k}(i)$.

```
X 7 8 6 1 7 8 6
8 X 5 7 8 2 5 7
6 8 X 5 6 8 3 5
7 5 6 X 7 5 6 4
5 3 4 2 X 3 4 2
4 6 1 3 4 X 1 3
2 4 7 1 2 4 X 1
3 1 2 8 3 1 2 X
```

Figure 5. The i,k entry is $\sigma_{8,k}(i)$.

```
 X 15 16 14  9 15 16 14  1 15 16 14  9 15 16 14
16  X 13 15 16 10 13 15 16  2 13 15 16 10 13 15
14 16  X 13 14 16 11 13 14 16  3 13 14 16 11 13
15 13 14  X 15 13 14 12 15 13 14  4 15 13 14 12
13 11 12 10  X 11 12 10 13 11 12 10  5 11 12 10
12 14  9 11 12  X  9 11 12 14  9 11 12  6  9 11
10 12 15  9 10 12  X  9 10 12 15  9 10 12  7  9
11  9 10 16 11  9 10  X 11  9 10 16 11  9 10  8
 9  7  8  6  1  7  8  6  X  7  8  6  1  7  8  6
 8 10  5  7  8  2  5  7  8  X  5  7  8  2  5  7
 6  8 11  5  6  8  3  5  6  8  X  5  6  8  3  5
 7  5  6 12  7  5  6  4  7  5  6  X  7  5  6  4
 5  3  4  2 13  3  4  2  5  3  4  2  X  3  4  2
 4  6  1  3  4 14  1  3  4  6  1  3  4  X  1  3
 2  4  7  1  2  4 15  1  2  4  7  1  2  4  X  1
 3  1  2  8  3  1  2 16  3  1  2  8  3  1  2  X
```

Figure 6. The i,k entry is $\sigma_{16,k}(i)$.